\documentclass[12pt]{article}
\usepackage{amssymb}
\usepackage{amssymb}
\usepackage{amsmath, amscd}
\usepackage{latexsym}
\usepackage[dvips]{graphicx}
\newtheorem{theorem}{Theorem}[section]
\newtheorem{lemma}[theorem]{Lemma}
\newtheorem{corollary}[theorem]{Corollary}

\newtheorem{assumption}[theorem]{Assumption}
\newtheorem{proposition}[theorem]{Proposition}

\newtheorem{question}{Question}[section]

\newcommand{\qed}{\hfill $\Box$ }
\newcommand{\proof}{\noindent{\bf Proof.}\ \ }
\begin{document}

\title{\Large {\bf Peripheral convex expansions of resonance graphs}}

\maketitle

\begin{center} {Zhongyuan Che \\ 
Department of Mathematics, Penn State University,\\
Beaver Campus,  Monaca, PA 15061, U.S.A.\\}
\end{center}

\begin{abstract}
In this paper, we show that the resonance graph
of a plane elementary bipartite graph $G$  can be obtained from an edge by a sequence of peripheral convex expansions 
with respect to a reducible face decomposition of $G$
if and only if the infinite face of $G$ is forcing.

\vskip 0.2in
\noindent {\emph{keywords}}: 
(peripheral) convex expansion, distributive lattice, forcing face, 
plane elementary bipartite graph, 
resonance graph, reducible face, $Z$-transformation graph
\end{abstract}

\section{Introduction}

Assume that $G$ is a connected graph and $H$ is a subgraph of $G$.
For any two vertices $x$ and $y$ of $H$,  if $d_H(x,y)=d_G(x,y)$, 
where $d_H(x,y)$ (resp., $d_G(x,y)$) is the length of a shortest path between $x$ and $y$ in $H$ (resp., in $G$),
then $H$ is \textit{isometric} in $G$;
if all shortest paths between $x$ and $y$ are contained in $H$, then $H$ is \textit{convex} in $G$.
Let $G_1$ and $G_2$ be isometric subgraphs of a connected graph $G$ such that 
$G_1 \cap G_2$ is a convex subgraph of $G$. 
Assume that the vertex set $V(G)=V(G_1) \cup V(G_2)$ and
there are no edges between $V(G_1) \setminus V(G_1 \cap G_2)$ and $V(G_2)\setminus V(G_1 \cap G_2)$.
Take disjoint copies of $G_1$ and $G_2$  
and connect every vertex of $G_1 \cap G_2$ in $G_1$ 
with the same vertex of $G_1 \cap G_2$ in $G_2$ with an edge. 
The we obtain a graph which is called the \textit{convex expansion} of $G$ with respect to $V(G_1)$ and $V(G_2)$.
In particular, if $G_1=G$ and  $G_2$ is convex in $G$,
then the above convex expansion  is called a \textit{peripheral convex expansion} of $G$.

A \textit{perfect matching} of a graph is a set of pairwise independent edges that covers all vertices of the graph.  
An edge of a graph is called \textit{allowed} if it is contained in a perfect matching, \textit{forbidden} otherwise.
A graph is \textit{elementary} if its allowed edges form a connected subgraph. 
The \textit{symmetric difference} of two sets is the set of elements belonging to the union but not the intersection of two sets.
The \textit{resonance graph} of a plane bipartite graph $G$, denoted by $Z(G)$,
is the graph whose vertices are the perfect matchings of $G$, 
and $M_1M_2$ is an edge of $Z(G)$ if the symmetric difference of $M_1$ and $M_2$ forms the boundary of a finite face $s$ of $G$,
and we say that $M_1M_2$ has the \textit{face label} $s$.
The concept of a resonance graph was introduced by many researchers independently under various names such that 
the $Z$-transformation graph and the perfect matching graph
because of its applications in chemistry, mathematics and physics.  The readers are referred to the survey paper \cite{Z06}.
A face of a plane bipartite graph is called a \textit{forcing face} if the boundary of the face is an even cycle 
and the subgraph obtained by removing all vertices of the face is either empty or has exactly one perfect matching \cite{CC13}. 

Let $G$ be a plane elementary bipartite graph and $s$ be a \textit{peripheral face} of $G$, 
that is, $s$ has some edges on the boundary of $G$. 
Let $H$ be the subgraph of $G$ obtained by removing all internal vertices (if exist) 
and edges on the common boundary of $s$ and $G$. 
If $H$ is elementary, then $s$ is called a \textit{reducible face} of $G$ \cite{TV12}. 
To provide approaches to an open question on characterizing $Z$-transformation graphs raised in the survey paper \cite{Z06}, 
we gave a convex expansion structure \cite{C18} of the resonance graph $Z(G)$ with respect to a reducible face of $G$.
It is well known \cite{LP86} that any elementary bipartite graph with more than two vertices is 2-connected.
We further showed \cite{C19} that if $G$ is a 2-connected outerplane bipartite graph, 
then $Z(G)$ can be obtained from an edge by a sequence of peripheral convex expansions 
with respect to a reducible face decomposition of $G$.

Motivated from the above work, we are interested in characterizing all resonance graphs which can be 
constructed from an edge by a sequence peripheral convex expansions.

\section{Preliminaries}\label{S:Preliminaries}

All graphs considered in this paper has a perfect matching unless specified otherwise.
Let $M$ be a perfect matching of a graph. 
An \textit{$M$-alternating cycle} (resp., $M$-alternating path)  is a cycle (resp., path) 
of the graph whose edges are alternately in and off $M$. 
An $M$-alternating cycle of a graph is also called a \textit{nice cycle}.
A path $P$ of the graph is called \textit{weakly $M$-augmenting} if either
$P$ is a single edge and not contained in $M$, 
or $P$ is an $M$-alternating path such that its two end edges are not contained in $M$. 
A face (including the infinite face) of a plane bipartite graph is called \textit{$M$-resonant} 
if its boundary is an $M$-alternating cycle for a perfect matching $M$ of the graph.
A plane bipartite graph $G$ is called \textit{weakly elementary} 
if and only if for each nice cycle $C$ of $G$, the subgraph of $G$ consisting of $C$ together with its interior is elementary.
Every plane elementary bipartite graph is weakly elementary \cite{ZZ00}. 

A \textit{reducible face decomposition} $G_i (1 \le i \le n)$ of 
a plane bipartite graph $G(=G_n)$ associated with the face sequence $s_i (1 \le i \le n)$
is defined as follows. Start from $G_1$ which is a finite face $s_1$ of $G$. If $G_{i-1}$ has been constructed, then
$G_i$ can be obtained from $G_{i-1}$ by adding an odd length path $P_i$  in the exterior of $G_{i-1}$ such that
only two end vertices of $P_i$ are contained in $G_{i-1}$,  
and $P_i$ and a part of the boundary of $G_{i-1}$ surround a finite face
$s_i$ of $G$ for all $2 \le i \le n$.  
By definition, we can see that for any reducible face decomposition $G_i (1 \le i \le n)$ of $G(=G_n)$ associated with 
the face sequence $s_i (1 \le i \le n)$, each $s_i$ is a reducible face of $G_i$ for $2 \le i \le n$. 
A plane bipartite graph with more than two vertices is elementary 
if and only if it has a reducible face decomposition \cite{ZZ00}. 

Let $S$ be a \textit{poset} with a partial order $\le$.  Assume that $x$ and $y$ are two elements of $S$.
If $x \le y$ and $x \neq y$, then we write $x<y$. 
We say that $y$ \textit{covers} $x$ if $x<y$ and there is no element $z \in S$ such that $x<z<y$.
The \textit{Hasse diagram} of $S$ is a digraph with the vertex set $S$ such that 
$\overrightarrow{yx}$ is a directed edge from $y$ to $x$ if and only if $y$ covers $x$.  
A \textit{lattice} is a poset with two operations $\vee$ and $\wedge$ 
such that $x \vee y$ is the unique least upper bound of $x$ and $y$,
and  $x \wedge y$ is the unique greatest lower bound of $x$ and $y$.
A lattice is \textit{distributive} if its two operations $\vee$ and $\wedge$  admit distributive laws.

A  \textit{median} of three vertices $x_1$, $x_2$ and $x_3$ of a graph $G$ is a vertex that is contained 
$I_{G}(x_1,x_2) \cap I_{G}(x_2,x_3) \cap I_{G}(x_1,x_3)$, where $I_{G}(x_i,x_j)$ is the set of all vertices on shortest
paths between $x_i$ and $x_j$ in $G$  for $1 \le i \neq j \le 3$. 
A graph is a \textit{median graph} if every triple vertices has a unique median.
It is well known  \cite{HIK11} that a graph is a median graph if and only if 
it can be obtained from a vertex by a convex expansion procedure. 

Assume that $G$ is a plane bipartite graph whose vertices are colored black and white 
such that adjacent vertices receive different colors. Let $M$ be a perfect matching of $G$.
An $M$-alternating cycle $C$ of $G$ is \textit{$M$-proper} (resp., \textit{$M$-improper}) 
if every edge of $C$ belonging to $M$ goes from 
white to black vertices (resp., from black to white vertices) along the clockwise orientation of $C$.
A plane bipartite graph $G$ with a perfect matching has a unique perfect matching $M_{\hat{0}}$ (resp., $M_{\hat{1}}$)
such that $G$ has no proper $M_{\hat{0}}$-alternating cycles (resp., no improper $M_{\hat{1}}$-alternating cycles) \cite{ZZ97}.
The \textit{$Z$-transformation digraph} (or, \textit{resonance digraph}), denoted by $\overrightarrow{Z}(G)$,
is the digraph obtained from $Z(G)$ by adding a direction for each edge  so that
$\overrightarrow{M_1M_2}$ is an directed edge from $M_1$ to $M_2$ if their symmetric difference 
is proper $M_1$-alternating (or, improper $M_2$-alternating).
Let $\mathcal{M}(G)$ be the set of all perfect matchings of $G$.
Then a partial order $\le$ can be defined on the set $\mathcal{M}(G)$ such that $M' \le M$ 
if there is a directed path from $M$ to $M'$ in $\overrightarrow{Z}(G)$.

Let $G$ be a plane (weakly) elementary bipartite graph. Then its resonance graph $Z(G)$ is a median graph \cite{ZLS08}, 
and $(\mathcal{M}(G), \le)$ is a finite distributive lattice
with the minimum $M_{\hat{0}}$ and the maximum $M_{\hat{1}}$, 
and the Hasse diagram of $(\mathcal{M}(G), \le)$ is isomorphic to $\overrightarrow{Z}(G)$  \cite{LZ03}.  
The \textit{height} of $(\mathcal{M}(G), \le)$ is defined as the length of a directed path from $M_{\hat{1}}$ to $M_{\hat{0}}$ in $\overrightarrow{Z}(G)$, 
we write $height(\mathcal{M}(G))$ briefly.
By the Jordan-Dedekind theorem \cite{B73} in finite distributive lattice,  all directed paths (if exist) between any pair of vertices 
of its Hasse diagram have the same length and are shortest paths between them. So, $height(\mathcal{M}(G))$ is 
the length of a shortest path between $M_{\hat{1}}$ and $M_{\hat{0}}$ in $Z(G)$ and can be written as
$height(\mathcal{M}(G)) = d_{Z(G)}(M_{\widehat{1}}, M_{\widehat{0}})$.

Let $G$ be a plane elementary bipartite graph and $\mathcal{M}(G)$ be the set of all perfect matchings of $G$. 
Assume that $s$ is a reducible face of  $G$.  By \cite{TV12},  the common boundary of $s$ and $G$ is an odd length path $P$. 
By Proposition 4.1 in \cite{C18},  $P$ is $M$-alternating for any perfect matching $M$ of $G$,
and $\mathcal{M}(G)=\mathcal{M}(G; P^{-}) \cup\mathcal{M}(G; P^{+})$,
where  $\mathcal{M}(G; P^{-})$ (resp., $\mathcal{M}(G; P^{+})$) is the set of perfect matchings $M$ of $G$ 
such that $P$ is weakly $M$-augmenting (resp., $P$ is not weakly $M$-augmenting).
Furthermore, $\mathcal{M}(G; P^{-})$ and $\mathcal{M}(G; P^{+})$ can be partitioned as  
\begin{eqnarray*}
\mathcal{M}(G; P^{-})&=&\mathcal{M}(G; P^{-}, \partial s) \cup\mathcal{M}(G; P^{-}, \overline{\partial s})\\
\mathcal{M}(G; P^{+})&=&\mathcal{M}(G; P^{+}, \partial s) \cup\mathcal{M}(G; P^{+}, \overline{\partial s})
\end{eqnarray*}
where $\mathcal{M}(G; P^{-}, \partial s)$ (resp., $\mathcal{M}(G; P^{-}, \overline{\partial s})$) 
is the set of perfect matchings $M$ in $\mathcal{M}(G; P^{-})$ such that $s$ is $M$-resonant (resp., not $M$-resonant),
and $\mathcal{M}(G; P^{+}, \partial s)$ (resp., $\mathcal{M}(G; P^{+}, \overline{\partial s})$) 
is the set of perfect matchings $M$ in $\mathcal{M}(G; P^{+})$ such that $s$ is $M$-resonant (resp., not $M$-resonant).

Two edges $uv$ and $xy$ of a graph $G$ are said to be in \textit{relation $\Theta$},
denoted by $uv \Theta xy$,  if $d_G(u,x)+d_G(v,y) \neq d_G(u,y)+d_G(v,x)$.
Relation $\Theta$ divides the edge set of a median graph into $\Theta$-classes 
as an equivalent relation \cite{HIK11}.
Based on the above partitions of $\mathcal{M}(G)$, 
a convex expansion structure of the resonance graph $Z(G)$ was given in \cite{C18}
using the tool of $\Theta$-relation.

\begin{assumption}\label{A:A}
Assume that $G$ is a plane elementary bipartite graph and $Z(G)$ is the resonance graph of $G$.
Let $s$ be a reducible face of $G$ and $P$ be the common boundary of $s$ and $G$.
Let $F$ be the set of all edges in $Z(G)$ with the face-label $s$.
Let $H$ be the subgraph of $G$ obtained by removing all internal vertices (if exist) and edges of $P$.
\end{assumption}

Let $\langle W \rangle$ represent the induced subgraph of a graph $G$ on $W \subseteq V(G)$.

\begin{theorem}\label{T:MedianZ(G)}\cite{C18}
Let Assumption \ref{A:A} hold true.
Then $F$ is a $\Theta$-class of $Z(G)$ and $Z(G)-F$ has exactly two components
$\langle \mathcal{M}(G; P^{-}) \rangle$ and $\langle \mathcal{M}(G; P^{+}) \rangle$.
Furthermore, (i) $F$ is a matching defining an isomorphism between 
$\langle \mathcal{M}(G; P^{-}, \partial s) \rangle$ and $\langle \mathcal{M}(G; P^{+}, \partial s) \rangle$;
 (ii) $\langle \mathcal{M}(G; P^{-}, \partial s) \rangle$ is convex in $\langle \mathcal{M}(G; P^{-}) \rangle$, 
$\langle \mathcal{M}(G; P^{+}, \partial s) \rangle$ is convex in $\langle \mathcal{M}(G; P^{+}) \rangle$;
(iii) $\langle \mathcal{M}(G; P^{-}) \rangle$ and $\langle \mathcal{M}(G; P^{+}) \rangle$ are median graphs, 
where $\langle \mathcal{M}(G; P^{-}) \rangle \cong Z(H)$.

In particular, $Z(G)$ can be obtained from $Z(H)$ by a peripheral convex expansion 
if and only if $\mathcal{M}(G; P^{+})=\mathcal{M}(G; P^{+}, \partial s)$.
\end{theorem}

\begin{proposition} \label{P:Sublattices}
Let Assumption \ref{A:A}  hold true.  
Then $(\mathcal{M}(G; P^{-}), \le)$, $(\mathcal{M}(G; P^{+}), \le)$, $(\mathcal{M}(G; P^{-}, \partial s), \le)$ 
and $(\mathcal{M}(G; P^{+}, \partial s), \le)$ 
are distributive sublattices of $(\mathcal{M}(G), \le)$.
\end{proposition}

\proof  By \cite{LZ03}, for any plane (weakly) elementary bipartite graph,
the set of its all perfect matchings forms a finite distributive lattice with
the Hasse diagram isomorphic to its resonance digraph. 
 
Since $H$ is a plane elementary bipartite graph, $\mathcal{M}(H)$ forms a finite distributive lattice 
with the Hasse diagram isomorphic to the resonance digraph $\overrightarrow{Z}(H)$.
By Theorem \ref{T:MedianZ(G)}, $\langle \mathcal{M}(G; P^{-}) \rangle \cong Z(H) = \langle \mathcal{M}(H) \rangle$.
Therefore, $(\mathcal{M}(G; P^{-}), \le)$ is  a distributive sublattice of $(\mathcal{M}(G), \le)$.

Remove two end vertices of $P$ and their incident edges.
If there is a pendent edge, then remove both end vertices of the pendent edge and their incident edges. 
Continue this way, if all vertices of $G$ can be removed, then $\langle \mathcal{M}(G; P^{+}) \rangle$ is the single vertex graph.
Otherwise, we obtain a nontrivial plane bipartite graph $H' \subset G$ without pendent edges.
Then $ \langle \mathcal{M}(H') \rangle = Z(H') \cong \langle \mathcal{M}(G; P^{+}) \rangle$. 
By Theorem \ref{T:MedianZ(G)}, $\langle \mathcal{M}(G; P^{+}) \rangle$ is a median graph and so is connected.
It follows that $Z(H')$ is connected.
By \cite{F03, ZZY04}, the resonance graph of a  plane bipartite graph  is connected 
if and only if the plane bipartite graph is weakly elementary. 
Hence, $H'$ is a plane weakly elementary bipartite graph. 
It follows that $\mathcal{M}(H')$ forms  a finite distributive lattice 
with the Hasse diagram isomorphic to the resonance digraph $\overrightarrow{Z}(H')$.
Therefore, $(\mathcal{M}(G; P^{+}), \le)$ is a distributive sublattice of $(\mathcal{M}(G), \le)$.

By Theorem \ref{T:MedianZ(G)}, $\langle \mathcal{M}(G; P^{-}, \partial s) \rangle$ is convex in $\langle \mathcal{M}(G; P^{-}) \rangle$, 
and $\langle \mathcal{M}(G; P^{+}, \partial s) \rangle$ is convex in $\langle \mathcal{M}(G; P^{+}) \rangle$, where 
$\langle \mathcal{M}(G; P^{-}) \rangle$ and $\langle \mathcal{M}(G; P^{+}) \rangle$ are median graphs.
By definition, $\langle \mathcal{M}(G; P^{-}, \partial s) \rangle$ (resp., $\langle \mathcal{M}(G; P^{+}, \partial s) \rangle$) 
is a median graph and so is a connected graph.
Similarly to the argument in the above paragraph, 
we can show that 
$\langle \mathcal{M}(G; P^{-}, \partial s) \rangle$ (resp., $\langle \mathcal{M}(G; P^{+}, \partial s) \rangle$)  is isomorphic to 
the resonance graph of a  plane weakly elementary bipartite graph $H_1 \subset G$ (resp., $H_2 \subset G$).
Therefore,  $(\mathcal{M}(G; P^{-}, \partial s), \le)$ and 
$(\mathcal{M}(G; P^{+}, \partial s), \le)$ are distributive sublattices of $(\mathcal{M}(G), \le)$.
\qed\\

\section{Main Results}

\begin{lemma}\label{L:MinMax-Z(G)-Z(H)}
Let Assumption \ref{A:A} hold true. Assume that $M_{\widehat{0}}$ and $M_{\widehat{1}}$ are the minimum 
and maximum of the finite distributive lattice $(\mathcal{M}(G), \le)$ on the set of all perfect matchings of $G$ respectively.
Then $Z(G)$ can be obtained from $Z(H)$ by a peripheral convex expansion
if and only if exactly one of $M_{\widehat{0}}$ and $M_{\widehat{1}}$ is contained in $\mathcal{M}(G; P^{+}, \partial s)$.
\end{lemma} 
\proof By definition of $\overrightarrow{Z}(G)$, all edges in $F$ are directed from one set  
of $\mathcal{M}(G; P^{-}, \partial s)$ and $\mathcal{M}(G; P^{+}, \partial s)$ to the other in $\overrightarrow{Z}(G)$. 
We distinguish two cases based on the above observation of edge directions.

Case 1. All edges in $F$ are directed from $\mathcal{M}(G; P^{-}, \partial s)$ to 
$\mathcal{M}(G; P^{+}, \partial s)$ in $\overrightarrow{Z}(G)$.
By Theorem \ref{T:MedianZ(G)},   $Z(G)-F$ has exactly two components
$\langle \mathcal{M}(G; P^{-}) \rangle$ and $\langle \mathcal{M}(G; P^{+}) \rangle$.
Then $M_{\widehat{1}} \in  \mathcal{M}(G; P^{-})$ and $M_{\widehat{0}} \in  \mathcal{M}(G; P^{+})$. 

 {\bf Claim. }  $\mathcal{M}(G; P^{+}, \partial s) =  \mathcal{M}(G; P^{+})$ 
if and only if $M_{\widehat{0}} \in  \mathcal{M}(G; P^{+}, \partial s)$.

Proof of Claim. Necessity is trivial.   
Sufficiency. Assume that $M_{\widehat{0}} \in \mathcal{M}(G; P^{+}, \partial s)$. 
Suppose that $\mathcal{M}(G; P^{+}, \partial s) \neq  \mathcal{M}(G; P^{+})$.
Then there is an $x \in \mathcal{M}(G; P^{+}) \setminus \mathcal{M}(G; P^{+}, \partial s)$.
Let $\overrightarrow{Q}$ be a directed path from $M_{\widehat{1}}$ passing through $x$ to $M_{\widehat{0}}$  in $\overrightarrow{Z}(G)$.
Let $Q$ be the corresponding path of $\overrightarrow{Q}$ in $Z(G)$.
By Theorem \ref{T:MedianZ(G)}, $Q$ contains an edge $M'_1M'_2 \in F$ satisfying
$M'_1 \in \mathcal{M}(G; P^{-}, \partial s)$ and $M'_2 \in \mathcal{M}(G; P^{+}, \partial s)$.
We can write \[\overrightarrow{Q}=\overrightarrow{Q}_{[M_{\widehat{1}}, M'_1]} \cup \overrightarrow{M'_1M'_2} \cup \overrightarrow{Q}_{[M'_2, x, M_{\widehat{0}}]},\]
where $\overrightarrow{Q}_{[M_{\widehat{1}}, M'_1]}$ (resp., $\overrightarrow{Q}_{[M'_2, x, M_{\widehat{0}}]}$)
is a directed path from $M_{\widehat{1}}$ to $M'_1$ (resp., from $M'_2$ passing through $x$ to $M_{\widehat{0}}$).
By Proposition \ref{P:Sublattices}, $(\mathcal{M}(G; P^{+}), \le)$ is a distributive sublattice of $(\mathcal{M}(G), \le)$.
By the Jordan-Dedekind theorem \cite{B73} in finite distributive lattice,  
all directed paths (if exist) between any pair of vertices 
of its Hasse diagram have the same length and are shortest paths between them.
Let $Q_{[M'_2, x, M_{\widehat{0}}]}$ be the corresponding path of $\overrightarrow{Q}_{[M'_2, x, M_{\widehat{0}}]}$ in $Z(G)$.
Then  $Q_{[M'_2, x, M_{\widehat{0}}]} \subseteq \langle \mathcal{M}(G; P^{+}) \rangle$
is a shortest path between  $M'_2$  and $M_{\widehat{0}}$ in $\langle \mathcal{M}(G; P^{+}) \rangle$. 
By Theorem \ref{T:MedianZ(G)}, $\langle \mathcal{M}(G; P^{+}, \partial s) \rangle$
is convex in the median graph $\langle \mathcal{M}(G; P^{+}) \rangle$,
any shortest path  between $M'_2 \in \mathcal{M}(G; P^{+}, \partial s)$  and 
$M_{\widehat{0}} \in \mathcal{M}(G; P^{+}, \partial s)$ is contained in $\langle \mathcal{M}(G; P^{+}, \partial s) \rangle$.
Hence, $Q_{[M'_2, x, M_{\widehat{0}}]} \subseteq \langle \mathcal{M}(G; P^{+}, \partial s) \rangle$.
This is a contradiction since $Q_{[M'_2, x, M_{\widehat{0}}]}$ contains a vertex 
$x \in \mathcal{M}(G; P^{+}) \setminus \mathcal{M}(G; P^{+}, \partial s)$. 
Therefore, $\mathcal{M}(G; P^{+}, \partial s) =  \mathcal{M}(G; P^{+})$.
This ends the proof of Claim.

By Theorem \ref{T:MedianZ(G)}, $Z(G)$ can be obtained from $Z(H)$ 
by a peripheral convex expansion with respect to $s$
if and only if $\mathcal{M}(G; P^{+}, \partial s) =  \mathcal{M}(G; P^{+})$
if and only if $M_{\widehat{0}} \in  \mathcal{M}(G; P^{+}, \partial s)$.

Case 2. All edges of $F$ are directed from $\mathcal{M}(G; P^{+}, \partial s)$ to 
$\mathcal{M}(G; P^{-}, \partial s)$ in $\overrightarrow{Z}(G)$.
Then $M_{\widehat{0}} \in  \mathcal{M}(G; P^{-})$ and $M_{\widehat{1}} \in  \mathcal{M}(G; P^{+})$. 
Similarly to the Claim in Case 1, we can show that 
$\mathcal{M}(G; P^{+}, \partial s) =  \mathcal{M}(G; P^{+})$ if and only if 
$M_{\widehat{1}} \in \mathcal{M}(G; P^{+}, \partial s)$.

By Theorem \ref{T:MedianZ(G)}, $Z(G)$ can be obtained from $Z(H)$ by a peripheral convex expansion with respect to $s$
if and only if $\mathcal{M}(G; P^{+}, \partial s) =  \mathcal{M}(G; P^{+})$
if and only if  $M_{\widehat{1}} \in \mathcal{M}(G; P^{+}, \partial s)$.
\qed\\

\begin{lemma} \label{L:Height-M(G)-M(H)}
Let Assumption \ref{A:A} hold true.  
Assume that $(\mathcal{M}(G), \le)$ and $(\mathcal{M}(H), \le)$ are the finite distributive lattices 
on the set of all perfect matchings of $G$ and $H$ respectively. 
Then $Z(G)$ can be obtained from $Z(H)$ by a peripheral convex expansion 
if and only if $height(\mathcal{M}(G))=height(\mathcal{M}(H)) + 1$.
\end{lemma}
\proof Necessity.  Assume that $Z(G)$ can be obtained from $Z(H)$ by a peripheral convex expansion with respect to $s$.
By Lemma \ref{L:MinMax-Z(G)-Z(H)}, without loss of generality, 
we can assume that $M_{\widehat{0}} \in \mathcal{M}(G; P^{+}, \partial s)$.
Let $M_0=M_{\widehat{0}} \oplus \partial s$ and $M_0|_H$ be the restriction of $M_0$ on $H$. 
Then $M_0|_H$ is a perfect matching of $H$.
Note that all $M_0|_H$-alternating cycles of $H$ are $M_0|_H$-improper.
By \cite{LZ03}, $M_0|_H$ is the minimum element of $(\mathcal{M}(H), \le)$.
It is easy to see that $M_0|_H$ can be extended uniquely to the perfect matching $M_0$ of $G$ in $\mathcal{M}(G; P^{-})$.  
By Proposition \ref{P:Sublattices}, $(\mathcal{M}(G; P^{-}), \le)$ is a distributive sublattice of $(\mathcal{M}(G), \le)$.
By Theorem \ref{T:MedianZ(G)}, $\langle \mathcal{M}(G; P^{-}) \rangle \cong Z(H)= \langle \mathcal{M}(H) \rangle$.
Then $M_0$ is the minimum element of the sublattice $(\mathcal{M}(G; P^{-}), \le)$.
By the construction of $M_0$, we can see that $s$ is proper $M_0$-resonant and $M_0 \oplus M_{\widehat{0}}=\partial s$.
Then $\overrightarrow{M_0M_{\widehat{0}}}$ is a directed edge in $\overrightarrow{Z}(G)$. 
Hence, there is a directed path $\overrightarrow{Q}$ from
$M_{\widehat{1}}$ to $M_{\widehat{0}}$ in $\overrightarrow{Z}(G)$ containing $\overrightarrow{M_0M_{\widehat{0}}}$. 
We can write $\overrightarrow{Q}=\overrightarrow{Q}_{[M_{\widehat{1}}, M_0]} \cup \overrightarrow{M_0M_{\widehat{0}}}$
where $\overrightarrow{Q}_{[M_{\widehat{1}}, M_0]}$ is a directed path from  $M_{\widehat{1}}$ to $M_0$.
Hence, \[height(\mathcal {M}(G)) = |\overrightarrow{Q}|
=|\overrightarrow{Q}_{[M_{\widehat{1}}, M_0]}|+ |\overrightarrow{M_0M_{\widehat{0}}}| 
= |\overrightarrow{Q}_{[M_{\widehat{1}}, M_0]}| +1.\]
Note that $|\overrightarrow{Q}_{[M_{\widehat{1}}, M_0]}|=height(\mathcal{M}(G; P^{-}))$
and $height(\mathcal{M}(G; P^{-}))=height(\mathcal{M}(H))$.
It follows that $height(\mathcal{M}(G))=height(\mathcal{M}(H)) + 1$. 

Sufficiency. 
Suppose that $Z(G)$ cannot be obtained from $Z(H)$ by a peripheral convex expansion with respect to $s$.
By Lemma \ref{L:MinMax-Z(G)-Z(H)},  
neither $M_{\widehat{0}}$ nor  $M_{\widehat{1}}$ is contained in $\mathcal{M}(G; P^{+}, \partial s)$.
Without loss of generality, we can assume that all edges in $F$
are directed from $\mathcal{M}(G; P^{-}, \partial s)$ to $\mathcal{M}(G; P^{+}, \partial s)$ in $\overrightarrow{Z}(G)$. 
Then $M_{\widehat{1}} \in \mathcal{M}(G; P^{-})$ and 
$M_{\widehat{0}} \in \mathcal{M}(G; P^{+}) \setminus \mathcal{M}(G; P^{+}, \partial s)$. 
Let $\overrightarrow{Q}$ be a directed path from $M_{\widehat{1}}$ to $M_{\widehat{0}}$ in $\overrightarrow{Z}(G)$. 
Let $Q$ be the corresponding path of $\overrightarrow{Q}$ in $Z(G)$.
By Theorem \ref{T:MedianZ(G)}, $Q$ contains an edge $xy \in F$ where $x \in \mathcal{M}(G; P^{-}, \partial s)$ and 
$y \in \mathcal{M}(G; P^{+}, \partial s)$. 
In particular, we can choose $x$ as the minimum element of the sublattice $(\mathcal{M}(G; P^{-}), \le)$.
Then $\overrightarrow{Q}=\overrightarrow{Q}_{[M_{\widehat{1}}, x]} \cup \overrightarrow{xy} \cup \overrightarrow{Q}_{[y, M_{\widehat{0}}]}$
where $\overrightarrow{Q}_{[M_{\widehat{1}}, x]}$ is a directed path 
from $M_{\widehat{1}}$ to $x$, 
and $\overrightarrow{Q}_{[y, M_{\widehat{0}}]}$ is a directed path 
from $y$ to $M_{\widehat{0}}$.
Note that $y \in \mathcal{M}(G; P^{+}, \partial s)$ is different from 
$M_{\widehat{0}} \in \mathcal{M}(G; P^{+}) \setminus \mathcal{M}(G; P^{+}, \partial s)$.
Then \[height(\mathcal{M}(G))= |\overrightarrow{Q}|
=|\overrightarrow{Q}_{[M_{\widehat{1}}, x]}|+ |\overrightarrow{xy}|+ |\overrightarrow{Q}_{[y, M_{\widehat{0}}]}| 
> |\overrightarrow{Q}_{[M_{\widehat{1}}, x]}| +1.\] 
Note that $|\overrightarrow{Q}_{[M_{\widehat{1}}, x]}|
=height(\mathcal{M}(G; P^{-}))=height(\mathcal{M}(H))$.
It follows that $height(\mathcal{M}(G))> height(\mathcal{M}(H))+1$. 
This is a contradiction to the assumption that $height(\mathcal{M}(G))=height(\mathcal{M}(H)) + 1$. 
\qed\\

Let $M_1$ and $M_2$ be two perfect matchings of a graph $G$.
Then a cycle of $G$ is called \textit{$(M_1, M_2)$-alternating} if its edges are in $M_1$ and $M_2$ alternately.
It is well know that the symmetric difference of $M_1$ and $M_2$, 
denoted by $M_1 \oplus M_2$,  is a union of vertex disjoint $(M_1,M_2)$-alternating cycles of $G$ \cite{LP86}.

\begin{figure}[h]
\begin{center}
\includegraphics[width=10cm]{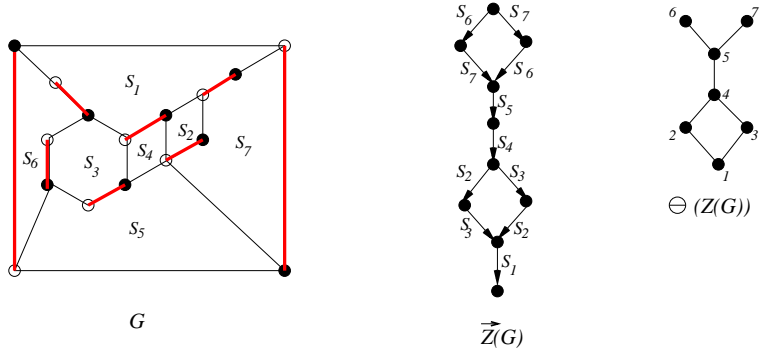}
\caption{An Example for Theorem \ref{T:PeripheralConvexExpansion}.} 
\label{F:Main-Example}
\end{center}
\end{figure}

\begin{theorem} \label{T:PeripheralConvexExpansion}  
Let $G$ be a plane elementary bipartite graph and $Z(G)$ be its resonance graph.
Then $Z(G)$ can be obtained from an edge by a sequence of peripheral convex expansions
with respect to a reducible face decomposition of $G$ if and only if the infinite face of $G$ is forcing.
\end{theorem}

\proof  Let $(\mathcal{M}(G), \le)$ be the finite distributive lattice on the set of all perfect matchings of $G$
with the maximum $M_{\widehat{1}}$ and the minimum $M_{\widehat{0}}$.
Let $\partial G$ denote the boundary of $G$.
By \cite{TV12}, $\partial G$ is both $M_{\widehat{0}}$-alternating and $M_{\widehat{1}}$-alternating.
Then $\partial G$ is an $(M_{\widehat{1}}, M_{\widehat{0}})$-alternating cycle. 
Let $\mathcal{F}$ be the set of all finite faces of $G$ and $f \in \mathcal{F}$.
Define $\phi_M(f)$ as the number of $(M, M_{\widehat{0}})$-alternating cycles 
in $M \oplus M_{\widehat{0}}$ containing $f$ in their interiors  \cite{ZLS08}.
Then $\phi_{M_{\widehat{0}}} (f)=0$ and $\phi_{M_{\widehat{1}}}(f)>0$ 
since $\partial G$ is an $(M_{\widehat{1}}, M_{\widehat{0}})$-alternating cycle.
By Theorem 3.2 in \cite{ZLS08} states that for any two vertices $M,M'$ of $Z(G)$,
distance $d_{Z(G)}(M, M') =\sum_{f \in \mathcal{F}} |\phi_{M} (f)- \phi_{M'} (f)|$.
Hence, $height(\mathcal{M}(G)) = d_{Z(G)}(M_{\widehat{1}}, 
M_{\widehat{0}}) =\sum_{f \in \mathcal{F}} |\phi_{M_{\widehat{1}}} (f)- \phi_{M_{\widehat{0}}} (f)|
=\sum_{f \in \mathcal{F}} \phi_{M_{\widehat{1}}} (f)$.
Therefore, we have the following height formula for the finite distributive lattice $(\mathcal{M}(G), \le)$.
\[height(\mathcal{M}(G)) =\sum_{f \in \mathcal{F}} \phi_{M_{\widehat{1}}} (f) \ \ \ \ \ \ \ \ \ \ (1) \]

Necessity. Assume that $Z(G)$ can be obtained from an edge by a sequence of peripheral convex expansions
with respect to a reducible face decomposition $G_i(1 \le i \le n)$ of $G(=G_n)$ where $n=|\mathcal{F}|$ 
is the number of finite faces of $G$.
By applying Lemma \ref{L:Height-M(G)-M(H)} repeatedly, we have 
$height(\mathcal{M}(G))=height(\mathcal{M}(G_{n-1}))+1=\cdots =height(\mathcal{M}(G_1))+(n-1)$.
Note that $height(\mathcal{M}(G_1))=1$ since $G_1$ is a finite face of $G$. 
Then $height(\mathcal{M}(G))=n$.
By formula (1), we have $height(\mathcal{M}(G))=n = \sum_{f \in \mathcal{F}} \phi_{M_{\widehat{1}}} (f)$
and so $\phi_{M_{\widehat{1}}}(f)=1$  for any $f \in \mathcal{F}$. 
Suppose that the infinite face of $G$ is not a forcing face.
Then the subgraph obtained by removing all vertices on $\partial G$ has at least two perfect matchings $M_1$ and $M_2$.
Recall that $M_1 \oplus M_2$ is a union of vertex-disjoint $(M_1, M_2)$-alternating cycles  \cite{LP86}.
Choose an $(M_1, M_2)$-alternating cycle $C$ such that $C$ is not contained in any other $(M_1, M_2)$-alternating cycles
in $M_1 \oplus M_2$.
Without loss of generality, we can assume that $C$ is proper $M_1$-alternating and so improper $M_2$-alternating.
It is clear that $M_1$ (resp., $M_2$) can be extended to a perfect matching $M'_1$ (resp., $M'_2$) of $G$ 
such that both $C$ and $\partial G$ are proper $M'_1$-alternating and improper $M'_2$-alternating,
and $C$ is not contained in any $(M'_1,M'_2)$-alternating cycles other than $\partial G$ in $M'_1 \oplus M'_2$.
Let $f \in \mathcal{F}$. Define $\psi_{M'_1M'_2} (f)$ as the number of proper $M'_1$-alternating cycles in $M'_1 \oplus M'_2$
with $f$ in their interiors minus the number of improper $M'_1$-alternating cycles in 
$M'_1 \oplus M'_2$ with $f$ in their interiors \cite{ZLS08}. 
It follows that there exists an $s \in \mathcal{F}$ contained in $C$ such that $\psi_{M'_1M'_2}(s) \ge 2$.
By Lemma 2.4 in  \cite{ZLS08}, $\phi_{M'_1}(s) - \phi_{M'_2}(s) = \psi_{M'_1M'_2}(s)$.
It follows that $\phi_{M'_1}(s) \ge 2$.  
Lemma 2.5 in  \cite{ZLS08} states that $M$ covers $M'$ in $\mathcal{M}(G)$ if and only if $\phi_{M}(f) - \phi_{M'}(f) = 1$
for $f=f_0$, where $f_0$ is a finite face bounded by the cycle $M \oplus M'$, and $0$ for the other faces $f$ in $\mathcal{F}$.
This implies that $\phi_{M_{\widehat{1}}}(s) \ge \phi_{M'_1}(s) \ge 2$ since 
$M_{\widehat{1}}$ is the maximum element of $\mathcal{M}(G)$.
On the other hand, we have shown that $\phi_{M_{\widehat{1}}}(s)=1$. 
This is a contradiction. Therefore, the infinite face of $G$ is forcing.

Sufficiency. If the infinite face of $G$ is forcing, then $M_{\widehat{0}} \oplus M_{\widehat{1}}=\partial G$
since  $\partial G$ is an $(M_{\widehat{1}}, M_{\widehat{0}})$-alternating cycle and 
$M_{\widehat{0}} \oplus M_{\widehat{1}}$ is a union of vertex-disjoint $(M_{\widehat{0}}, M_{\widehat{1}})$-alternating cycles.
By definition, $ \phi_{M_{\widehat{1}}} (f)=1$ for each $f \in \mathcal{F}$.
By formula (1), we have $height(\mathcal{M}(G)) = \sum_{f \in \mathcal{F}} \phi_{M_{\widehat{1}}} (f) =\sum_{f \in \mathcal{F}} 1=|\mathcal{F}|=n$.
Let $G_i (1 \le i \le n)$ be a reducible face decomposition of $G$ associated with the face sequence $s_i (1 \le i \le n)$. 
By Lemma \ref{L:Height-M(G)-M(H)}, for $2 \le i \le n$,
$height(\mathcal{M}(G_i)) \ge height(\mathcal{M}(G_{i-1}))+1$ and the equality holds if and only if 
$Z(G_i)$ can be obtained from $Z(G_{i-1})$ by a peripheral convex expansion.
Recall that $height(\mathcal{M}(G_1))=1$. 
Then $height (\mathcal{M}(G))=n$ if and only if $Z(G)$ can obtained from an edge by a sequence of peripheral convex expansions
with respect to a reducible face decomposition of $G$.
\qed\\

\begin{corollary}\label{C:PCE}
Let $G$ be a plane elementary bipartite graph such that each finite face has a vertex on the boundary of $G$.
Then $Z(G)$ can be obtained from an  edge by a sequence of peripheral convex expansions
with respect to a reducible face decomposition of $G$.
\end{corollary}

For example, the plane elementary bipartite graph $G$ given in Figure 23 \cite{CT14}
obtained from a six crossing knot universe  
satisfies the property of Corollary \ref{C:PCE}.  
Its resonance graph $Z(G)$ can be obtained from an edge by a sequence of peripheral convex expansions 
with respect to a reducible face decomposition of $G$ associated with the face sequence $s_i (1 \le i \le 7)$.

\begin{figure}[h]
\begin{center}
\includegraphics[width=9cm]{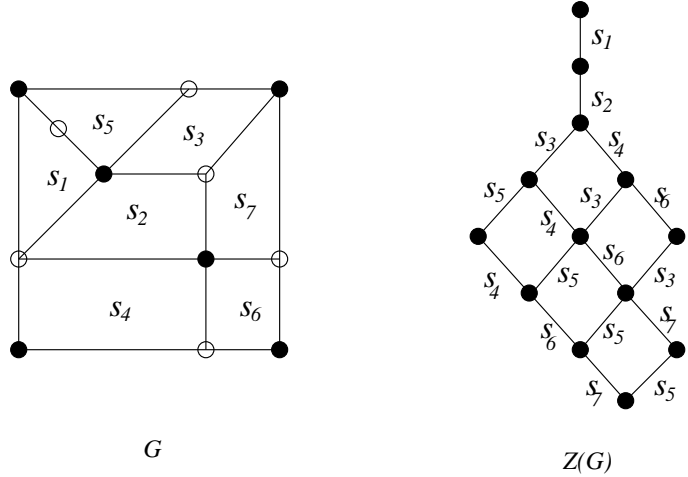}
\caption{An example for Corollary \ref{C:PCE} from Figure 23 in \cite{CT14}.} 
\label{F:Corollary-Example}
\end{center}
\end{figure}

The \textit{induced graph $\Theta(Z(G))$} on the $\Theta$-classes of the edge set of a resonance graph $Z(G)$ is the graph  
whose vertex set is the set of $\Theta$-classes, and two $\Theta$-classes $E_1$ and $E_2$ are adjacent if  
$Z(G)$ has two incident edges $e_1 \in E_1$ and $e_2 \in E_2$ such that $e_1$ and $e_2$  
are not contained in a common 4-cycle of $Z(G)$ \cite{V05}. 

In  \cite{C19}, we showed that if $G$ is a 2-connected outerplane bipartite graph, then $Z(G)$ can be obtained from an edge by a sequence of peripheral convex expansions
with respect to a reducible face decomposition of $G$. Moreover,  $\Theta(Z(G))$ is a tree and isomorphic to the inner dual of $G$, 
This generalized the corresponding result in \cite{V05} for catacondensed hexagonal systems. 
Note that for a plane elementary bipartite graph $G$, if $Z(G)$ can be obtained from an edge by a sequence of peripheral convex expansions
with respect to a reducible face decomposition of $G$, then it is not necessarily true that $\Theta(Z(G))$ is isomorphic to the inner dual of $G$.
See Examples given in Figure \ref{F:Main-Example} and Figure \ref {F:Corollary-Example}. 
\begin{question}
When a plane elementary bipartite graph $G$ has the property that $\Theta(Z(G))$ is isomorphic to the inner dual of $G$?
\end{question}

  

\end{document}